\documentclass[12 pt]{amsart}

\usepackage{amscd, amssymb}

\newtheorem{pr}{Proposition}

\newtheorem{tm}{Theorem}

\newcommand{\proj}{\mathbf P}

\newcommand{\rarr}{\rightarrow}
\newcommand{\oh}{{\mathcal{O}}}
\newcommand{\com}{\mathbb{C}}
\newcommand{\Q}{\mathbb{Q}}
\newcommand{\Z}{\mathbb{Z}}

\newcommand{\eqq}{\stackrel{\sim}{=}}

\def\scap{\mathbin{\text{\scriptsize$\cap$}}}

\newcommand{\bpf}{\noindent {\em Proof.} }
\newcommand{\epf}{\qed \vspace{+10pt}}

\begin{document}
\title{Hodge integrals and degenerate contributions}
\author{R. Pandharipande}
\date{25 March  1999}
\maketitle

\pagestyle{plain}
\setcounter{section}{-1}
\section{\bf{Introduction}}
\subsection{}

Let $X$ be a  nonsingular, projective, 3 dimensional complex algebraic
variety.
Let $\overline{M}_{g_D,n}(X,\beta)$ be the moduli space of
stable maps from genus $g_D$ curves to $X$ representing
the homology class $\beta \in H_2(X, \Z)$.
The Gromov-Witten invariants of $X$ are defined via tautological
integrals
over these moduli spaces of maps (against their virtual
fundamental classes):
$$N^{g_D}_\beta(\gamma_1, \ldots, \gamma_n)
 = \int_{[\overline{M}_{g_D,n}(X,\beta)]^{vir}} \prod_{i=1}^n
\text{ev}_i^*(\gamma_i),$$
where $\text{ev}_i$ is the $i^{th}$ evaluation map and
$\gamma_i \in H^*(X,\Z)$.
As the moduli spaces are Deligne-Mumford stacks, the
Gromov-Witten invariants take values in $\Q$.
Let $T_X$ and $K_X$ be the tangent bundle and the canonical class
of $X$.
For a 3-fold, the dimension
formula shows the virtual dimensions 
do not depend upon the genus:
$$\text{dim}_{{vir}}(\overline{M}_{g_D}(X,\beta))
= 3g_D-3+ \chi(T_X)=-K_X \cdot \beta.$$
If we restrict attention to a fixed curve 
class $\beta \in H_2(X,\Z)$, there are two basic
possibilities: $-K_X\cdot \beta =0$ or $-K_X\cdot \beta>0$
(the negative case is of no interest here since then
the Gromov-Witten invariants
vanish). We will always take $\beta \neq 0$.

\subsection{Case $-K_X \cdot \beta=0$.}
If $X$ is Calabi-Yau, this case holds for all classes $\beta$.
Let $d$ be a positive integer.
Let $C\subset X$ be a nonsingular genus $g< g_D$ curve of
class $\beta/d$. 
The moduli space $\overline{M}_{g_D}(X,\beta)$ contains a 
substack of maps with genus $g_D$ domains which factor through a $d$-fold
cover of $C$. Under suitable conditions, this substack of maps
covering  $C$ is a connected 
component of $\overline{M}_{g_D}(X,\beta)$. 
In the latter case, the contribution
of $C$ to the genus $g_D$, class $\beta$ 
Gromov-Witten invariant of $X$ is well-defined. It is these 
degenerate contributions
that are studied here.
Degenerate contributions play a central role in identifying
the integer quantities in the Gromov-Witten theory of $X$. These
integrality properties remain a very mysterious part of the
subject.
In algebraic geometry, degenerate contributions are related to
Hodge integrals over the moduli space of curves $\overline{M}_{g,n}$ [FP].
In string theory, recent progress in the study of  these
contributions has been made by a link to M-theory [GV1], [GV2] (see also
[MM]). While the mathematical results presented  here overlap 
with the M-theoretic results of [GV2], the precise
connection between the two approaches is still not completely
understood. The differences are discussed below in Section \ref{mmmm}.

Let $C\subset X$ be a nonsingular genus $g$ curve representing
the class $[C]\in H_2(X, \Z)$. For the degenerate analysis,
we assume the normal bundle to $C$ in $X$ is general.
Consider the moduli space of maps 
$\overline{M}_{g+h}(X,d[C])$.
If $g=0$ or $1$,  this moduli space
will have a connected component equal to $\overline{M}_{g+h}(C, d[C])$.
The contribution $C_g(h,d)$ of $C$
to the genus $g+h$ Gromov-Witten is thus 
well-defined for $g=0,1$ and
all values $h\geq 0$, $d>0$. 
The above component claim relies on rigidity arguments which
possibly fail for multiple covers of genus $g\geq 2$ curves.
However, in the degree 1 case, $\overline{M}_{g+h}(X,[C])$
has a component equal to $\overline{M}_{g+h}(C,[C])$
for all $g$ and $h$. Hence, $C_g(h,1)$ is always
well-defined. At present, because of the possibility of deformations
in $X$ away from $C$,
the correct definition of $C_g(h,d)$ in general 
is not known to the author.

The contributions in case $g=0$ have recently been
calculated in algebraic geometry [FP] and string 
theory [GV1], [MM]:
\begin{equation}
\label{klkl}
\sum_{h=0}^{\infty} C_0(h,1) t^{2h} = 
\Big( \frac{\sin(t/2)}{t/2}  \Big)^{-2},
\end{equation}
\begin{equation}
\label{fsfs}
C_0(h,d) = d^{2h-3} C_0(h,1),
\end{equation}
where $C\subset X$ is a nonsingular, rigid rational curve. 
The contribution $$C_0(0,d)=1/d^3$$ is the Aspinwall-Morrison
formula which had been proven previously by several
different methods [AM], [M], [V]. 

A nonsingular curve $C\subset X$ is {\em rigid} if $H^0(C, N)=0$ where
$N$ is the normal bundle of $C$ in $X$.
For rational $C$, rigidity is 
equivalent to the bundle splitting $N = \oh(-1) \oplus \oh(-1)$.
Define
$C\subset X$ to be {\em super-rigid} if, for
all non-constant maps of nonsingular curves $\mu: C' \rarr C$,
$$H^0(C', \mu^*(N))=0.$$
It is clear rigidity and super-rigidity are equivalent 
in the rational case, but differ for higher genus.
Super-rigidity is a generic condition on the
normal bundle for elliptic
curves in $X$.
Kley has informed the author his existence
result for rigid elliptic curves on complete intersection Calabi-Yau
3-folds 
 also shows the existence of super-rigid
elliptic curves [K].

The contributions $C_1(0,d)$ are easily computed for
super-rigid elliptic curves $C$.
The component of the moduli space
$\overline{M}_1(X,d[C])$ corresponding to 
maps with image $C$ is nonsingular
of dimension 0 (and equal to $\overline{M}_1(C,d[C])$ ).
The points of $\overline{M}_1(C, d[C])$
correspond naturally to the set of subgroups of
$\Z \oplus \Z$ of index $d$.  Hence, after accounting
for automorphisms, 
$$C_1(0,d) = \frac{\sigma(d)}{d} = \sum_{i|d} \frac{1}{i}$$
(see, for example, [S]).
All other contributions of an elliptic curve $C$
vanish by the following result.

\begin{tm}
Let $C\subset X$ be a super-rigid elliptic curve. Then,
$$C_1(h,d)=0$$ 
for all $h>0$, $d>0$.
\end{tm}
\noindent This vanishing
was conjectured by Gopakumar-Vafa in [GV1] and
is derived in M-theory in [GV2]. The proof given here
uses basic constructions related to the virtual
fundamental class.

The degree 1 contributions $C_g(h,1)$ take a very simple form.
\begin{tm} Let $g\geq 0$.
$$\sum_{h=0}^{\infty} C_g(h,1) t^{2h} = 
\Big( \frac{\sin(t/2)}{t/2}  \Big)^{2g-2}.$$
\end{tm}
\noindent Theorem 2 is derived in Section 2 by expressing the
contributions $C_g(h,1)$ as Hodge integrals
over the moduli space of curves.
The required integrals are then computed via geometric
constructions, relations, and series manipulations.
Theorem 2 is the main result of this paper.

The right side of Theorem 2 was encountered before in
the following related result of [FP]:
\begin{equation}
\label{www}
1+ \sum_{h\geq 1} \sum_{i=0}^{h} t^{2h} k^i
\int_{\overline{M}_{h,1}} \psi_1^{2h-2+i} \lambda_{h-i} =
\Big( \frac{\sin(t/2)}{t/2} \Big)^{-k-1}.
\end{equation}
Theorem 2 gives an interpretation of these
Hodge integrals in the Gromov-Witten theory of 
Calabi-Yau 3-folds.

\subsection{M-theory predictions}
\label{mmmm}
The method of [GV1], [GV2] is to consider limits of type IIA string
theory which may be conjecturally analyzed in M-theory.
A remarkable proposal is made in [GV2] for the form of the
Gromov-Witten potential $\tilde{F}$ of a Calabi-Yau 3-fold $X$.
Let $$\tilde{F}(t,q)= 
\sum_{g\geq 0} t^{2g-2} \tilde{F}_g(t,q),$$
$$\tilde{F}_g(t,q)= 
\sum_{0 \neq \beta \in H_2(X, \mathbb{Z})} N_{\beta}^g\ q^\beta,$$
where $N_\beta^g$ is the genus $g$ Gromov-Witten
invariant of $X$ in curve class $\beta$. The potential
$\tilde{F}$ differs from the {\em full} potential by the
constant map $(\beta=0)$ contribution -- the constant contributions
have been calculated in [FP], [GV1], [MM].
For each curve class $\beta\in H_2(X, \mathbb{Z})$ and
genus $g$, there is an {\em integer} $n_\beta^g$ counting BPS states
in the associated M-theory. The formula of [GV2] is:
\begin{equation}
\label{ttttt}
\tilde{F}(t,q)
= \sum_{g,\beta} n_\beta^g t^{2g-2} \sum_{d>0} \frac{1}{d}\Big(
\frac{\text{sin}
({dt/2})}{t/2}\Big)^{2g-2} q^{d\beta}.
\end{equation}

If $C^M_g(h,d)$ denotes the contribution of a single BPS state
in genus $g$ and class $\beta$ to the Gromov-Witten
invariant in genus $g+h$ and class $d\beta$, then formula
(\ref{ttttt}) is equivalent to the equations:
$$\sum_{h=0}^\infty C^M_g(h,1)t^{2h} = 
\Big( \frac{\sin(t/2)}{t/2}  \Big)^{2g-2},$$
$$C^M_g(h,d)= d^{2g+2h-3} C^M_g(h,1).$$
The first of these agrees with Theorem 2, so
$C^M_g(h,1)=C_g(h,1)$.
The second is a generalization
of (\ref{fsfs}) to $g\geq 0$. It is therefore reasonable to interpret the
states $n_\beta^0$ as counting embedded (virtual) curves of genus $0$ and
degree $\beta$ (even for the
Calabi-Yau quintic these numbers $n_\beta^0$ are at best virtual
because of the existence of Vainsencher's nodal rational curves).
However, when specialized to genus 1,
the second equation yields $C^M_1(0,d)= 1/d$ in
contrast to $C_1(0,d)= \sigma(d)/d$. The (virtual) count of embedded
genus 1 curves should be derived from $\tilde{F}_1$ via the
multiple cover corrections $C_0(1,d)$ and $C_1(0,d)$ 
(as previously pursued in [BCOV]).
Gromov-Witten theory would predict the resulting number to
be virtually enumerative, and thus integral (this heuristic
argument for integrality is not a proof). 
The M-theoretic perspective predicts a {\em different}
correction of $\tilde{F}_1$ to yield integers via formula
(\ref{ttttt}). Klemm has checked the two genus 1
integrality predictions
both hold in low degrees for several Calabi-Yau 3 folds [Kl]. 
These integrality constraints are not trivially dependent.
No proofs of any of these integrality constraints
are known to the author.

To find higher genus evidence for the formula (\ref{ttttt}),
a direct computation of the potential
$\tilde{F}$ in the local Calabi-Yau case ($\proj^2$ with canonical
bundle) for low genera and degrees has been pursued by
Klemm and Zaslow [KlZ].
The Gromov-Witten invariants (in all genera) may be computed
in this case by the virtual localization formula of [GP] and
the holomorphic anomally equation [BCOV]. The integrality
predicted by (\ref{ttttt}) is a nontrivial constraint which 
is verified in all calculations. 

At this point, it is not clear how to define or compute the
general contributions $C_{g}(h,d)$. One may hope a complete
understanding of $C_g(h,d)$ 
will lead to an integrality property of the
Gromov-Witten potential of $X$ distinct from (\ref{ttttt}).

\subsection{Case $-K_X \cdot \beta>0$}
In this case, the moduli spaces $\overline{M}_{g_D}(X, \beta)$
have positive virtual dimensions. The Gromov-Witten invariants
$N_\beta^{g_D}(\gamma)$ of $X$ then depend
upon a vector of cohomology classes 
$$\gamma=(\gamma_1, \ldots, \gamma_k), \ \ 
\gamma_i \in H^*(X,\Z).$$
Let $Y_i \subset X$ be general topological
cycles dual to the classes $\gamma_i$.
If we wish to identify integers in this Gromov-Witten
theory, degenerate contributions again  play a role.

Let us assume we are in an ideal situation with respect
to the moduli spaces of maps to $X$.
Let $M^{Bir}_{g}(X, \beta)$ denote the
moduli space of birational maps from smooth genus $g$ domain
curves.   
We
assume first:
\begin{enumerate}
\item[(i)]
The spaces $M^{Bir}_{g}(X,\beta)$ 
are generically reduced and of the expected dimension for
all $g\leq g_D$.
\end{enumerate}
There is then an enumerative integer
$E_\beta^{g_D}(\gamma)$ defined
to equal the number of genus $g_D$ maps of class
$\beta$ with smooth domains meeting all the cycles $Y_i$.
However, $E_\beta^{g_D}(\gamma)$ will not equal $N_\beta^{g_D}(\gamma)$.
The difference arises from the following observation.
For each $g< g_D$, there are $E_\beta^{g}(\gamma)$ maps with 
smooth genus $g$ domains of class $\beta$ satisfying the required
incidence conditions. The Gromov-Witten invariant
$N_\beta^{g_D}(\gamma)$ receives a degenerate
contribution from each of these lower genus solutions (via
reducible
genus $g_D$ maps factoring through covers of the lower genus curves).
As the genus $g$ solution represents the class $\beta$,
the covers must be of degree 1. 
These degenerate contributions are therefore analogous to
$C_g(g_D-g,1)$.

Dimension counts show maps multiple onto
their image and maps 
with reducible images
are not {\em expected} to contribute to $N_\beta^{g_D}(\gamma)$.
This is the second ideal assumption:
\begin{enumerate}
\item[(ii)] Maps in $\overline{M}_{g_D}(X,\beta)$
multiple onto their image or
with reducible image do not satisfy incidence conditions
to the cycles $Y_i$.
\end{enumerate}

Let $C\subset X$ be a nonsingular, genus $g$ curve of
class $\beta$ satisfying incidence conditions
to the cycles $Y_i$. Assume further $C$ is infinitesimally rigid 
with respect to these incidence conditions.
The contribution
$C_g(h,X,\beta)$ of $C$
to the Gromov-Witten invariant
$N_\beta^{g+h}(\gamma)$ is then well-defined: it is
found in Section 3 to be an integral over
the moduli space $\overline{M}_{g+h}(C,[C])$.
This contribution is  easily seen to be independent of $\gamma$.
The final ideal assumption is:
\begin{enumerate}
\item[(iii)] For all $g<g_D$, the solution maps
counted by $E^g_\beta(\gamma)$ are nonsingular
embeddings infinitesimally rigid with respect
to the incidence conditions.
\end{enumerate}
The ideal relation between Gromov-Witten theory and
the enumerative invariants is:
\begin{equation}
\label{iidd}
N^{g_D}_\beta(\gamma)= 
\sum_{g=0}^{g_D} C_g(g_D-g, X, \beta) E_\beta^g(\gamma).
\end{equation}
The validity of the relation (\ref{iidd}) for  
$N^{g_D}_\beta(\gamma)$
requires assumptions (i), (ii), and (iii) for
$g_D$, $\beta$, and $\gamma$.

The easiest 3-fold to consider is $X=\proj^3$. As the divisor
$-K_{\proj^3}$ is ample, $-K_{\proj^3} \cdot \beta>0$
for all nonzero curve classes.
The moduli spaces of maps to $\proj^3$ are easily seen to
be ideal in the above sense for the genera $g_D=0,1,2$, all
degrees $d>0$, and all $\gamma$. 
The rigidity statements follow as usual from
Bertini arguments (see [FuP]).
Therefore, the ideal relation (\ref{iidd}) holds in these
genera. The equation
\begin{equation}
\label{fdfd} 
N_d^0= E_d^0
\end{equation}
is well known for $\proj^3$ (we drop 
$\gamma$ in these equations).
In joint work with Getzler and  Graber,
we had computed
$$C_0(1,\proj^3,d)=  \frac{1-2d}{12},$$
\begin{equation}
\label{gg11}
N_d^1= \frac{1-2d}{12} E^0_d + E^1_d.
\end{equation}
Equation (\ref{gg11}) was used in Getzler's study [Ge] of the
genus 1 enumerative geometry of $\proj^3$.
Using Xiong's calculations of low degree 
genus 2 Gromov-Witten invariants
of $\proj^3$ as data,  Jinzenji and  Xiong conjectured the contribution
equation:
\begin{equation}
\label{gg22}
{N}_{d}^{2}=\frac{3 -11d+10d^2}{720} E_d^{0}-
\frac{4d}{24} E_d^{1}+E_d^{2}.
\end{equation}

These equations led Jinzenji and  Xiong to recently conjecture a
general formula [J] analogous to Theorem 2:
\begin{equation}
\label{jinzx}
\sum_{h=0}^{\infty} C_g(h,X,\beta) t^{2h} = 
\Big( \frac{\sin(t/2)}{t/2}  \Big)^{2g-2-K_X\cdot \beta}.
\end{equation}
The contribution $C_g(h,X,\beta)$ is calculated here by the method
used in the proof of Theorem 2. 
\begin{tm} The degenerate contributions $C_g(h,X,\beta)$ are
determined by formula (\ref{jinzx}).
\end{tm}

Theorem 3 and relation (\ref{iidd}) prove
formulas (\ref{fdfd}), (\ref{gg11}), (\ref{gg22}) for $g=0,1,2$ and 
all degrees $d>0$ in $\proj^3$. 
For higher genera, it is known the space of curves in $\proj^3$
may be of excess dimension. For example, the moduli space
$\overline{M}_{3}(\proj^3,4)$ has a 17 dimensional
component, but is expected to be 16 dimensional.
The definition of
enumerative invariants is therefore not clear from a space curve point
of view.
However, the invariants
$E^g_\beta(\gamma)$ may still be {\em defined} by Theorem 3 from
the Gromov-Witten invariants and equation (\ref{iidd}).
Perhaps an integrality property
holds for $E^g_\beta(\gamma)$ in some general context.

Algebraic 3-folds are special in Gromov-Witten theory
since the (virtual) dimensions of the moduli spaces of stable maps 
do not depend upon the genus. A similar uniform treatment
of degenerate contributions
in higher dimensions will require new ideas. 
Graber has carried out a related degenerate analysis
in the genus 0 Gromov-Witten theory of the
Hilbert scheme of 2 points of $\proj^2$ [Gr].

\subsection{Moduli of curves}
The Hodge integral approach taken
here has an application to the geometry of the
moduli space of nonsingular curves $M_g$,  $(g\geq 2)$.
The tautological
ring $\mathcal{R}^*(M_g)$ is the subring of 
$A^*(M_g)$ generated by the $\kappa$ classes (see [Mu]).
The intersection calculus of $\mathcal{R}(M_g)$
has a very rich structure. A detailed study by  Faber
of $\mathcal{R}(M_g)$ for low genera
has led to very precise conjectures of this ring
structure [F1]. In particular, 
Faber 
has conjectured $\mathcal{R}^*(M_g)$ is a
Gorenstein ring with socle in degree $g-2$.
In [GeP], the (conjectural) intersection pairing of $\mathcal{R}(M_g)$
is directly linked to Gromov-Witten theory via  
(conjectural) Virasoro constraints on
the descendent potential of $\proj^2$. The computation here
of the degenerate contributions $C_g(h,1)$  leads to   
a formula in $\mathcal{R}^*(M_g)$
conjectured previously by Faber from evidence for 
$g\leq15$.

\begin{tm} For $g\geq 2$, the relation
$$ \sum_{i=0}^{g-2} (-1)^i \lambda_i \kappa_{g-2-i} = \frac{2^{g-1}}{g!} 
\kappa_{g-2}$$
holds in $\mathcal{R}^*(M_g)$.
\end{tm}

The author thanks P. Belorousski, 
C. Faber, E. Getzler, T. Graber, M. Jinzenji, S. Katz,
A. Klemm, H. Kley, C. Vafa, and E. Zaslow  for
comments and correspondence related to degenerate contributions.
In particular, this paper was inspired by questions of C. Vafa.
The author was partially support by 
National Science
Foundation grant DMS-9801574.

\section{\bf Theorem 1}
\subsection{Super-rigidity}
Let $C\subset X$ be a nonsingular elliptic
curve in a Calabi-Yau 3-fold. The normal bundle $N$
is of rank 2 with trivial determinant. If $C$ is rigid, a straightforward
argument shows 
$N$ contains a non-trivial degree 0 line sub-bundle
$L$:
$$0 \rarr L \rarr N \rarr L^{-1} \rarr 0.$$
Conversely, such a filtration implies the rigidity of $C$.
The curve $C$ is super-rigid
if and only if $L$ is not a torsion element of
the Picard group of $C$. 
While super-rigidity is
a stronger condition on $N$ than rigidity, it is an open 
condition.
Super-rigidity is required for the
equality of moduli spaces proven in Proposition 1.
Note super-rigidity implies $H^0(C', \mu^*(N))=0$ for
every non-constant {\em stable} map $\mu:C' \rarr C$. 

The moduli spaces $\overline{M}_{1+h}(C,d[C])$ and
$\overline{M}_{1+h}(X,d[C])$ are Deligne-Mumford stacks
with possibly nonreduced structures.

\begin{pr} 
Let $C\subset X$ be a nonsingular, super-rigid elliptic
curve. The space of maps 
$\overline{M}_{1+h}(C,d[C])$ is a union of connected components
of $\overline{M}_{1+h}(X, d[C])$ for all $h\geq 0$, $d>0$.
\end{pr}

\bpf There is a natural map:
$$\iota: 
\overline{M}_{1+h}(C,d[C]) \rarr \overline{M}_{1+h}(X, d[C]).$$
By the super-rigidity of $C$, the locus of
$\overline{M}_{1+h}(X, d[C])$ corresponding to
maps with support in $C$ is a union of connected
components of $\overline{M}_{1+h}(X, d[C])$. 
We will prove $\iota$ is an isomorphism onto these
connected components.

It suffices to prove a lifting property for families of
stable maps over Artinian local rings $A$. Let
$\xi \in {\text{Spec}}(A)$ be the geometric point corresponding
to the maximal ideal $m\subset A$. Let
$$\pi: \mathcal{F} \rarr {\text{Spec}}(A), \ \ \mu: \mathcal{F} \rarr X$$
be a family of stable maps
satisfying
\begin{equation}
\label{poq}
\mu_\xi: \mathcal{F}_\xi \rarr C \subset X.
\end{equation}
We will prove $\mu$ factors through $C$. This lifting
implies the desired isomorphism property of $\iota$.

Let $\mathcal{I}$ be the ideal sheaf of $C$ in $X$.
We must prove the natural map  
$$\phi:\mu^*(\mathcal{I}) \rarr \oh_{\mathcal{F}}$$
is zero.  Certainly
$\phi$ has image in $m \oh_{\mathcal{F}}$ by the
assumption (\ref{poq}) on the geometric fiber $\xi$.
Hence, $\phi$ induces a natural map on $\mathcal{F}$:
\begin{equation}
\label{qdq}
\mu^*(\mathcal{I}/ \mathcal{I}^2) \rarr 
m \oh_{\mathcal{F}} /m^2 \oh_{\mathcal{F}} =
(m/m^2) \otimes_\com \oh_{\mathcal{F}_\xi}.
\end{equation}
The restriction of $\mu^*(\mathcal{I}/ \mathcal{I}^2)$
to $\mathcal{F}_{\xi}$ is simply $\mu_\xi^*(N^*)$.
By the super-rigidity of $C$, the map (\ref{qdq}) is zero. 
We conclude
$\phi$ factors through $m^2 \oh_{\mathcal{F}}$.

The above argument may be used to prove the
following implication: if $\phi$ factors through
$m^k \oh_{\mathcal{F}}$, then $\phi$ factors through
$m^{k+1} \oh_{\mathcal{F}}$.
Since $A$ is Artinian, $m$ is nilpotent. Hence,
$\phi$ vanishes.
\epf

There are two perfect obstruction theories on $\overline{M}_{1+h}(C,
d[C])$ obtained from considering the moduli problem of
maps to $C$ and $X$ respectively (see [B], [BF], [LT]).
Let $$\pi: \mathcal{F} \rarr \overline{M}_{1+h}(C,d[C]),$$
$$\mu: \mathcal{F} \rarr C$$
be the universal family and universal map respectively.
By super-rigidity $\pi_*\mu^*(N)=0$ and
$R^1\pi_*\mu^*(N)$ is a rank $2h$ bundle.
The two obstruction theories differ exactly
by the bundle $R^1\pi_*\mu^*(N)$.
From the definition of the virtual class, we conclude:
\begin{equation}
\label{formm}
C_1(h,d)= \int_{[\overline{M}_{1+h}(C, d[C])]^{vir}}
c_{2h}(R^1\pi_*\mu^*(N)).
\end{equation}

\subsection{Vanishing results}
Let $E$ be any bundle on $C$.
Consider the complex $R\pi_*\mu^*(E)$ in the
derived category of coherent sheaves on $\overline{M}_{1+h}(C,d[C])$.
Let $\mathcal{L}$ be a $\pi$-relatively ample polarization on
$\mathcal{F}$. We may find an exact sequence of bundles on
$\mathcal{F}$:
$$0 \rarr K \rarr \oplus \mathcal{L}^{-k} \rarr \mu^*(E) \rarr 0$$
for some positive integer $k$ [H]. 
As both $\pi_*K$ and $\pi_* \mathcal{L}^{-k}$ vanish, we find a 
two term bundle resolution of $R\pi_*\mu^*(E)$:
$$[R^1\pi_* K \rarr R^1\pi_* \oplus \mathcal{L}^{-k}] \eqq
R\pi_*\mu^*(E).$$
The Chern classes of $R\pi_*\mu^*(E)$ are defined
by $c(R^1\pi_* K)/ c(R^1\pi_* \oplus \mathcal{L}^{-k})$.
This definition is independent of two term resolutions in
the derived category.

As $\pi_*\mu^*(N)=0$ and
$R^1\pi_*\mu^*(N)$ is a rank $2h$ bundle,
we see (\ref{formm}) may now be rewritten as:
\begin{equation*}
C_1(h,d)= \int_{[\overline{M}_{1+h}(C, d[C])]^{vir}}
[c^{-1}(R\pi_*\mu^*(N))]_{2h}.
\end{equation*}

It is easy to find flat families of bundles on $C$
connecting $N$ and the the trivial rank 2
bundle $I=\oh_C \oplus \oh_C$.
 For example, if $P$ is a sufficiently ample line bundle,
both $N\otimes P$ and $I\otimes P$ will have nowhere
vanishing sections:
$$0 \rarr \oh_C \rarr N \otimes P \rarr P^2 \rarr 0,$$
$$0 \rarr \oh_C \rarr I \otimes P \rarr P^2 \rarr 0.$$
Hence $N$ and $I$ are connected in the
family of extensions of $P$ by $P^{-1}$.
The integral
$$\int_{[\overline{M}_{1+h}(C, d[C])]^{vir}}
[c^{-1}(R\pi_*\mu^*(E))]_{2h}$$
is clearly constant as $E$ varies in this family (for
example, the two term resolutions of $R\pi_*\mu^*(E)$
may be chosen compatibly over the family).
We conclude,
$$C_1(h,d)= \int_{[\overline{M}_{1+h}(C, d[C])]^{vir}}
[c^{-1}(R\pi_*\mu^*(I))]_{2h}$$

Now assume $h>0$.
Let $\gamma: \overline{M}_{1+h}(C, d[C]) \rarr \overline{M}_{1+h}$
be the natural map to the moduli space of curves.
Let $\mathbb{E}$ denote the Hodge bundle on $\overline{M}_{1+h}$:
the fiber of $\mathbb{E}$ over the moduli point 
$[F]\in \overline{M}_{1+h}$ is $H^0(F, \omega_F)$ (see [Mu]).
Since 
$$\pi_*\mu^*(I) = \oh_{\overline{M}} \oplus \oh_{\overline{M}},$$
$$R^1\pi_* \mu^*(I)= \gamma^*(\mathbb{E}^* \oplus \mathbb{E}^*),$$
we see $[c^{-1}(R\pi_*\mu^*(I))]_{2h}$
is a cohomology class pulled-back via $\gamma$ from 
$\overline{M}_{1+h}$.
Hence, to complete the proof of
 Theorem 1, it suffices to show the following vanishing.
\begin{pr} 
Let $h>0$. Then, $\gamma_*([\overline{M}_{1+h}(C,d[C])]^{vir})=0$.
\end{pr}
\bpf
Fix a base point $p\in C$ for the course of the proof.
We will consider the moduli space of 1-pointed maps
$\overline{M}_{1+h,1}(C,d[C])$. Let
$$\text{ev}_1^{-1}(p)
=\overline{M}_{1+h,p}(C,d[C]) \subset \overline{M}_{1+h,1}(C, d[C])$$
denote the subspace of maps for which the marking has image $p$.
There is a canonical isomorphism obtained by the
group law on $C$:
$$ \overline{M}_{1+h,1}(C, d[C]) \eqq C \times
\overline{M}_{1+h,p}(C, d[C]).$$
Let $\rho: \overline{M}_{1+h,1}(C, d[C]) \rarr
\overline{M}_{1+h,p}(C,d[C])$ denote the
canonical projection.

The perfect obstruction theory on $\overline{M}_{1+h,1}(C, d[C])$
may be obtained from a canonical distinguished triangle involving
the cotangent complex of the Artin stack of prestable curves and
the perfect obstruction theory relative to this Artin stack
(see [B], [BF], [GrP]). These objects are naturally equivariant
with respect to the natural group law on $C$ (see the constructions
of [B], [BF]). Hence, the virtual class of $\overline{M}_{1+h,1}(C,d[C])$
is a pull-back of an algebraic cycle class on 
$\overline{M}_{1+h,p}(C, d[C])$. As the map
$$\gamma_1: \overline{M}_{1+h,1}(C, d[C]) \rarr \overline{M}_{1+h,1}$$
factors through $\overline{M}_{1+h,p}(C, d[C])$, we obtain
\begin{equation}
\label{qas}
\gamma_{1*}( [\overline{M}_{1+h,1}(C, d[C])]^{vir}) =0.
\end{equation}

Consider now the commutative diagram obtained from the
$1$-pointed moduli spaces:
\begin{equation}
\label{fred}
\begin{CD}
 \overline{M}_{1+h,1}(C,d[C])
@>>{\gamma_1}>  \overline{M}_{1+h,1} \\
@V{\pi}VV   @V{\pi}VV \\
 \overline{M}_{1+h}(C, d[C]) @>>{\gamma}> \overline{M}_{1+h}.
\end{CD}
\end{equation}
While (\ref{fred}) is not a fiber square, it is easy to
see the following equality holds
\begin{equation}
\label{asd}
\gamma_{1*} \pi^* = \pi^* \gamma_*.
\end{equation}
From the Axiom of contracting a point [BM], we see
$$\pi^*( [\overline{M}_{1+h}(C,d[C])]^{vir})= 
[\overline{M}_{1+h,1}(C, d[C])]^{vir}.$$
Then, equations (\ref{qas}) and (\ref{asd}) imply:
\begin{equation}
\label{xcv}
\pi^* \gamma_*([\overline{M}_{1+h}(C,d[C])]^{vir}) =0.
\end{equation}
For any class $\xi \in A_*(\overline{M}_{1+h})$,
$$\pi_*( \psi_1 \cdot \pi^*(\xi))= 2h \cdot \xi,$$
where $\psi_1$ is the Chern class of the cotangent line
on $\overline{M}_{1+h,1}$.
Hence, the pull-back
$\pi^*: A_*(\overline{M}_{1+h}) \rarr A_*(\overline{M}_{1+h,1})$
is injective. The Proposition now follows from
(\ref{xcv}).
\epf

\section{\bf Theorem 2}
\subsection{Rigidity}
Let $C\subset X$ be a rigid, nonsingular genus $g$ curve with
normal bundle $N$. 
The contribution $C_g(0,1)$ is certainly 1, so we may assume
$h$ is a positive integer. The proof of Proposition 1 also establishes the
following facts. First,
the moduli space $\overline{M}_{g+h}(C, [C])$ is a component
(easily seen to be connected) of
$\overline{M}_{g+h}(X, [C])$. Second, the contribution $C_g(h,1)$
is determined by:
\begin{equation}
\label{fortt}
C_g(h,1)= \int_{[\overline{M}_{g+h}(C, [C])]^{vir}}
c_{2h}(R_{g,h}).
\end{equation}
Here, $R_{g,h}$ denotes the rank $2h$ bundle $R^1\pi_*\mu^*(N)$.
Note the virtual dimension of
$\overline{M}_{g+h}(C, [C])$ is also $2h$.
The arguments of Section 1 are valid because a rigid curve
is super-rigid in degree 1.

\subsection{Irreducible components of $\overline{M}_{g+h}(C, [C])$}
Let $C$ be a nonsingular genus $g$ curve. 
Let $h$ be a positive integer.
We first analyze the
moduli space of degree 1 maps $\overline{M}_{g+h}(C, [C])$. 
Let $P(h)$ denote the set of partitions $h$.
There is a natural set-theoretic function:
$$\nu: \overline{M}_{g+h}(C, [C]) \rarr P(h)$$
defined by the following method.
Let $\mu: F \rarr C$ correspond to a point 
$[\mu]\in \overline{M}_{g+h}(C,[C])$. The domain $F$ must contain
a unique irreducible component $F_C$  mapped isomorphically to 
$C$ by $\mu$. The arithmetic genera of the
connected components $\{F_i\}$ of $\overline{F \setminus F_C}$
form  a partition of $h$. Let $\nu([\mu])$ equal this
partition. The irreducible components of $\overline{M}_{g+h}(C,[C])$
are in bijective correspondence with $P(h)$ by the
value of $\nu$ on a general element. 

Let $\tau=(h_1\geq \ldots \geq  h_l)$
 be a partition of $h$ of length $l$. 
Consider the Fulton-MacPherson configuration space $C[l]$ of
$l$ marked points in $C$:
$C[l]$ is a natural compactification
of the space of $l$ distinct points of $C$ [FuM]. 
If $C$ has no automorphisms,
$C[l]$ is simply the fiber of $\overline{M}_{g,l} \rarr \overline{M}_g$
over the moduli point $[C]$. Define the nonsingular Deligne-Mumford
stack $I_\tau$ by:
\begin{equation}
\label{zxz}
I_\tau = C[l] \times \prod_{i=1}^l \overline{M}_{h_i,1}.
\end{equation}
There is a natural family, 
$$\pi: \mathcal{F} \rarr I_\tau,$$
of prestable curves over
$I_\tau$ obtained by 
attaching a $1$-pointed
genus $h_i$ curve to the $i^{th}$ marking of 
the universal $l$-pointed curve over $C[l]$.
Moreover, there is canonical projection $\mu: \mathcal{F} \rarr C$.
The induced morphism:
$$ \gamma_\tau: I_\tau \rarr 
\overline{M}_\tau \subset \overline{M}_{g+h}(C, [C])$$
is finite and surjective onto the irreducible
component $\overline{M}_\tau$ corresponding to the partition $\tau$.

Let $\partial \overline{M}_{h_i,1}$ denote the boundary
of the moduli space: the locus of curves with at least one node.
Similarly, let $\partial C[l]\subset C[l]$ 
denote the locus of nodal curves
($\partial C[l]$ may also 
be viewed as the locus lying over the diagonals
of the product $C^l$). Let $\partial I_\tau$ denote
the 
union of the pull-backs of the boundaries of the factors
(\ref{zxz}) via the $l+1$
projections. Let 
$$\partial\gamma_\tau: \partial I_\tau \rarr \overline{M}_{g+h}(C,[C])$$
denote the natural map.
The main geometric result used in the proof
of Theorem 2 is the following vanishing.

\begin{pr} For all partitions $\tau$ of $h$,
$$c_{2h}(\partial \gamma_\tau ^*(R_{g,h})) =0.$$
\end{pr}

\bpf
By definition, $\partial I_\tau$ is union of the
pull-backs of the boundary divisors of the $l+1$ product factors
of (\ref{zxz}).
We show $c_{2h}(\partial\gamma_\tau^* (R_{g,h}))$ restricts
to 0 on each of these pull-backs.

Let $\text{pr}_j$ denote the projection of
$I_\tau$ onto the $(j+1)^{st}$ factor of (\ref{zxz})
for $0\leq j \leq l$.
There are $l$ natural evaluation maps $\text{ev}_i: C[l]\rarr C$
obtained from the $l$ markings.
Define $\text{ev}_i: I_\tau \rarr C$ 
by the composition
$$I_\tau \stackrel{\text{pr}_0}{\longrightarrow} 
C[l] \ \  \stackrel{\text{ev}_i}
\longrightarrow \ \ C$$
for $1\leq i \leq l$.

The bundle $\gamma_\tau^*(R_{g,h})$ is easily analysed
via   
the natural normalization sequence of the family
$\mathcal{F}$. We find a decomposition:
\begin{equation}
\label{mmast}
\gamma_\tau^*(R_{g,n}) = \bigoplus_{i=1}^l
\mathbb{E}_i^* \otimes \text{ev}_i^*(N)
\end{equation}
where $\mathbb{E}_i$ is the Hodge bundle on $\overline{M}_{h_i,1}$.
We denote the pull-back of these
Hodge bundles to $I_\tau$ by the same symbols.
An important relation among the Chern classes of the
Hodge bundle has been established by Mumford in [Mu].
Mumford's relation 
is:
$c(\mathbb{E}_i)\cdot c(\mathbb{E}_i^*)=1$ in $A^*(\overline{M}_{h_i,1}).$

From (\ref{mmast}), we deduce:
\begin{equation*}
c_{2h}(\gamma_\tau^*(R_{g,h}))=  \prod_{i=1}^{l} c_{2h_i}
(\mathbb{E}_i^* \otimes \text{ev}_i^*(N)).
\end{equation*}
Algebra and Mumford's relation
then yield:
\begin{equation}
\label{qqqq}
c_{2h}(\gamma_\tau^*(R_{g,h}))=  \prod_{i=1}^{l} \lambda_{h_i}
\lambda_{h_i-1} c_1(\text{ev}_i^*(N^*)).
\end{equation}
Here, $\lambda_k$ denotes the $k^{th}$ Chern class of
the Hodge bundle. 

First, consider a boundary 
divisor $\Delta \subset \overline{M}_{h_j,1}$.
The pull-back of $\Delta$ to $I_\tau$ is simply:
$$\text{pr}_j^{-1}(\Delta)= 
C[l] \times \Delta
\times \prod_{i\neq j} \overline{M}_{h_i,1}.$$
The restriction of the factor $\lambda_{h_j} \lambda_{h_j-1}$
of (\ref{qqqq}) to $\Delta$ has been proven by Faber
to vanish [F1]
(the reducible
divisors of $\overline{M}_{h_1,1}$ have non-trivial
genus splittings). Hence, the
restriction of $c_{2h}(\gamma_\tau^*(R_{g,h}))$ to $\text{pr}_j^{-1}
(\Delta)$ vanishes.

Second, consider a boundary divisor $\Delta$ of
$C[l]$. The divisor $\Delta$ corresponds to a locus
in which a subset $S\subset [l]$ (of at least 2 elements) of 
the marked points coincide over $C$.
The evaluation maps $\{ \text{ev}_i \}_{i\in S}$
coincide when restricted to $\text{pr}_{0}^{-1} (\Delta)$.
Therefore, since $c_1(N^*)^2=0$, the
restriction of $c_{2h}(\gamma_\tau^*(R_{g,h}))$ to $\text{pr}_j^{-1}
(\Delta)$ vanishes.
\epf

\subsection{Hodge integrals}
Let $\partial \overline{M}_\tau= \gamma_\tau(\partial I_\tau)$, and
let $M_\tau= \overline{M}_\tau \setminus \partial \overline{M}_\tau$. 
$M_\tau$ is open in $\overline{M}_{g+h}(C,[C])$ 
and corresponds to the moduli space of degree 1 maps which
consist of nonsingular curves of genus $h_i$ attached to
distinct point of $C$. 

A deformation theory argument
shows $M_\tau$ is a nonsingular moduli stack of
dimension $\sum_{i=1}^l (3h_i-1)$. 
More precisely, for $[\mu: F \rarr C]\in M_\tau$,
there is a canonical exact sequence:
$$0 \rarr \text{Aut}_{[F]} \rarr H^0(F, \mu^*(T_C)) \rarr \text{Def}_{[\mu]} 
\stackrel{\iota}{\rarr}$$
$$\ \ \ 
\text{Def}_{[F]} \rarr H^1(F, \mu^*(T_C) \rarr \text{Obs}_{[\mu]}\rarr 0$$
where $\text{Aut}_{[F]}$ is the infinitesimal automorphisms of
$F$ and
$\text{Def}_{[F]}$,
$\text{Def}_{[\mu]}$ are the infinitesimal deformation
spaces of $F$, $\mu$ respectively.
It is easy to prove the
cokernel of $\iota$ is equal to a vector space $V$ with filtration
$$0\rarr \text{Def}_{[C]} \rarr V \rarr 
\bigoplus_{i=1}^l (T_{p_i} \otimes T_{p'_i}) \rarr 0.$$
Here, the component $F_i \subset F$ of genus $h_i$ is attached to
$C$ at the points $p_i \in F_i$ and $p'_i \in C$.
The cokernel computation amounts to showing the map $\mu$ has no 
infinitesimal deformations which smooth
any of the $l$ nodes of $F$.
We then see $\text{Def}_{[\mu]}$ is of constant
dimension $\sum_{i=1}^l (3h_i-1)$. Moreover, the
obstruction space is a bundle over $M_\tau$ with fiber
\begin{equation}
\label{qty}
\frac{ H^1(F, \mu^*(T_C))}  {Im(V)} 
= \bigoplus_{i=1}^l \frac{ H^1(F_i, \oh_{F_i} 
\otimes T_{p'_i})}{ T_{p_i} \otimes T_{p'_i}}
\end{equation}
over $[\mu]$.
The essential point here is the
deformation theory of maps in $M_\tau$ is very
simple.

Let $\text{Aut}_\tau$ denote the stabilizer of the
permutation $\mathbb{S}_l$-action 
on the $l$-tuple $\tau$. 
The map $\gamma_\tau: I_\tau \rarr 
\overline{M}_\tau$ is $\text{Aut}_\tau$-invariant.
Moreover, the quotient map
induces a proper, bijective morphism
\begin{equation}
\label{rew}
\tilde{\gamma}_\tau: I_\tau/\text{Aut}_\tau \rarr \overline{M}_\tau.
\end{equation}
Let $I^0_\tau= I_\tau \setminus \partial I_\tau$.
Certainly, $\tilde{\gamma}_\tau$ induces an isomorphism
$I^0_\tau/\text{Aut}_\tau \eqq M_\tau$.

The restriction of 
the virtual class $\xi^{vir}=[\overline{M}_{g+h}(C,[C])]^{vir}$
to the disjoint open union $\bigcup_{\tau \in P(h)} M_\tau$
is:
$$ \bigoplus_{\tau \in P(h)} \xi^{vir}_\tau,$$
where $\xi^{vir}_\tau \in A_{2h}(M_\tau)$.
The pull-back of $\xi^{vir}_\tau$ to $I_\tau^0$ is identified
from the obstruction theory (\ref{qty}) to be:
\begin{equation}
\label{vbv}
\gamma_\tau ^*(\xi^{vir}_\tau) =
\prod_{i=1}^{l} c_{h_i-1} \Big( \frac{c(\mathbb{E}_i^* \otimes
\text{ev}_i^*(T_C))}{1-\psi_1+c_1(\text{ev}_i^*(T_C))}
\Big).
\end{equation}
Since $M_\tau$ is nonsingular, the restriction of the virtual class is
the Euler class of the
obstruction bundle.

The virtual class $\xi^{vir}$ may be (non-canonically) expressed
at a sum:
$$\bigoplus_{\tau \in P(h)} \overline{\xi}^{vir}_\tau,$$
where $\overline{\xi}^{vir}_\tau \in A_{2h}(\overline{M}_\tau)$.
Using the proper bijection (\ref{rew}),
we see:
\begin{equation}
\label{hgh}
C_g(h,1)= \sum_{\tau\in P(h)} \int_{I_\tau/\text{Aut}_\tau}
\overline{\xi}^{vir}_\tau \scap c_{2h}( \tilde{\gamma}_\tau^* (R_{g,h})).
\end{equation}
By the vanishing of Proposition 3, equation (\ref{hgh})
remains valid if $\overline{\xi}_\tau^{vir}$ is replaced with {\em any}
cycle class which restricts to $\xi_\tau^{vir}$ on $M_\tau$.
This observation together with (\ref{vbv}) yields the
equality:
\begin{equation}
\label{hghh}
C_g(h,1)= \sum_{\tau\in P(h)} \frac{1}{|\text{Aut}_\tau|} 
\int_{I_\tau}
c_{2h}({\gamma}_\tau^* (R_{g,h})) \cdot
\prod_{i=1}^l
  c_{h_i-1} \Big( \frac{c(\mathbb{E}_i^*\otimes \text{ev}_i^*(T_C))}
{1-\psi_1+c_1(\text{ev}_i^*(T_C))}
\Big).
\end{equation}
Equation (\ref{qqqq}) together with basic algebraic
manipulations then prove the main integral formula:
\begin{equation}
\label{hkl}
 C_g(h,1)= \sum_{\tau\in P(h)} \frac{(2-2g)^l}{|\text{Aut}_\tau|} 
\prod_{i=1}^l
\int_{\overline{M}_{h_i,1}} 
\lambda_{h_i} \lambda_{h_i-1} (\sum_{j=0}^{h_i-1}(-1)^j
\lambda_j \psi_1^{h_i-1-j}).
\end{equation}
The only aspect of $N$ which affects the integral
\begin{equation*}
C_g(h,1)= \int_{[\overline{M}_{g+h}(C, [C])]^{vir}}
c_{2h}(R_{g,h}).
\end{equation*}
is $\int_C c_1(N^*)$. This Chern class enters 
enters (\ref{hkl}) via equation
(\ref{qqqq}) yielding the factor $$(\int_C c_1(N^*))^l = (2-2g)^l.$$
Theorems 2-4 will directly follow from  formula (\ref{hkl}).

For $q\geq 1$, define $\alpha_q = \int_{\overline{M}_{q,1}}
\lambda_{q}\lambda_{q-1} (\sum_{j=0}^{q-1}(-1)^j
\lambda_j \psi_1^{q-1-j})$.
Define the generating series:
$$Q(t)=  \sum_{q\geq 1} \alpha_q t^{2q}.$$
An immediate consequence of formula (\ref{hkl}) is the equation:
\begin{eqnarray*}
\sum_{h\geq 0} C_g(h,1) t^{2h} &  =  & \text{exp}((2-2g)Q(t)) \\
& = & \text{exp}(2Q(t))^{1-g} \\
&= & (\sum_{h\geq 0} C_0(h,1) t^{2h})^{1-g} \\
& = & \Big( \frac{\sin(t/2)}{t/2}  \Big)^{2g-2}.
\end{eqnarray*}
The last equality follows from the previous computations of
$C_0(h,1)$ in [FP].
The proof of Theorem 2 is complete.

\section{\bf Theorem 3}
We follow here the notation of Section 0.4 .
Let
$C$ be a nonsingular genus $g$ curve in
a 3-fold $X$ representing the homology class $\beta$.
We now assume $-K_X\cdot \beta>0$, so the moduli space
$\overline{M}_{g}(X,\beta)$ is of positive expected
dimension. Let $\gamma=(\gamma_1, \ldots, \gamma_n)$
be a vector of cohomology classes defining a Gromov-Witten
invariant $N_\beta^g(\gamma)$.
For each $i$, Let $Y_i\subset X $ be a topological cycle dual to 
$\gamma_i$. 
Let $p_i \in C \cap Y_i$.
We let $(C)$  denote the identity map $\pi: C \rarr C \subset X$
defining a point in the moduli space of stable maps.
The contribution of $C$ to $N_{\beta}^{g+h}(\gamma)$ via
covers
   will
require two general position hypotheses analogous to
rigidity in the Calabi-Yau case:
\begin{enumerate}
\item[(i)] $(C, p_1,\ldots, p_n)$ is a nonsingular
            point of $\overline{M}_{g,n}(X,\beta)$ lying on a component of
             expected dimension $-K_X \cdot \beta+n$.
\item[(ii)] The topological intersection of the
            cycles $\text{ev}_i^{-1}(Y_i)$ in
            $\overline{M}_{g,n}(X,\beta)$ 
            is transverse at $(C,p_1, \ldots
               p_n)$.
\end{enumerate}
Under these hypotheses, the degenerate contribution of
$C$ may be expressed directly as an integral over
$\overline{M}_{g+h}(C,[C])$.

Let $W\subset \overline{M}_g(X,\beta)$ be the
open, nonsingular, expected dimensional subset of the
moduli space of maps.
Let $U \subset W$ be the open
subset  
corresponding to embeddings of
nonsingular genus $g$ curves in $X$. 
As such embeddings have no nontrivial automorphisms, $U$ is
a nonsingular variety (not just a Deligne-Mumford stack). Moreover,
by assumption (i), $U$ is 
nonempty of dimension $-K_X\cdot \beta$ and contains $(C)$.
After discarding a finite number of points of $U$, we may assume
$(C)$ is the only point of $U$ meeting all the cycles $Y_i$.
Note the moduli space $U$ is also naturally an open set
of a component of the Hilbert scheme of curves in $X$.
Let 
$$\eta:\mathcal{C} \rarr U$$
denote the universal family of curves over $U$.
Let $\overline{M}_{g+h}(\eta, \beta)$ denote the
$\eta$-relative moduli space of maps representing the
fundamental class of the fibers of $\eta$. 
There is a natural morphism of Deligne-Mumford stacks:
\begin{equation}
\label{qrtqrt}
\iota: \overline{M}_{g+h}(\eta, \beta) \rarr \overline{M}_{g+h}
(X,\beta)
\end{equation}
obtained by composition. There are several tautological
morphisms (over $U$):
$$ \pi: \mathcal{F} \rarr \overline{M}_{g+h}(\eta,\beta),$$
$$ \mu: \mathcal{F} \rarr \mathcal{C},$$
$$ \tau: \overline{M}_{g+h}(\eta,\beta) \rarr U.$$
Let $\mathcal{N}$ denote the
universal normal bundle $\mathcal{N}$ on $\mathcal{C}$.
$\mathcal{N}$ is the family of normal bundles of the fibers
of $\eta$ in $X$. As $U$ is nonsingular of expected dimension,
$\eta_*(\mathcal{N})$ is isomorphic to the tangent
bundle of $U$ and $R^1\eta_*(\mathcal{N})=0$.

A deformation theoretic
check over Artinian rings shows $\iota$ is an open immersion.
We see the stack $\overline{M}_{g+h}(\eta,\beta)$ has
two natural fundamental classes. The first is 
$[\overline{M}_{g+h}(\eta,\beta)]^{vir}$
obtained from the structure
of a $\eta$-relative moduli space of maps. Second, the open 
inclusion $\iota$ endows $\overline{M}_{g+h}(\eta, \beta)$
with the perfect obstruction theory on $\overline{M}_{g+h}(X,\beta)$.
A direct comparison of these two obstruction theories on
$\overline{M}_{g+h}(\eta, \beta)$ shows they differ
exactly by the bundle $R_{g,h}=R^1\pi_*\mu^*(\mathcal{N})$:
\begin{equation}
\label{jkljkl}
\iota^* ([\overline{M}_{g+h}(X,\beta)]^{vir}) =
[\overline{M}_{g+h}(\eta,\beta)]^{vir} \scap c_{2h}(R_{g,h}).
\end{equation}
Relations
 (\ref{qrtqrt}) and (\ref{jkljkl}) are valid when considered
in the context of $n$-pointed stable maps (this may be deduced
from the above unpointed relations together with the natural properties
of these virtual structures under the morphisms forgetting the
markings [BM]).

By relation (\ref{jkljkl}) and the
definition of the Gromov-Witten invariants,
the
contribution of $(C, p_1,\ldots, p_n)$ to $N^{g+h}_\beta(\gamma)$
is equal to the intersection product:
\begin{equation}
\label{sqsqsq}
[\overline{M}_{g+h,n}(\eta,\beta)]^{vir} \scap c_{2h}(R_{g,h})\scap
\prod_{i=1}^n \text{ev}_i^{-1}(Y_i),
\end{equation}
with value in the zeroth homology of the compact space
$\cap_{i=1}^n \text{ev}_i^{-1}(Y_i).$
By assumption (ii) and
the pull-back properties of the virtual class, intersection
(\ref{sqsqsq}) is (numerically) equal to:
\begin{equation}
\label{sqsqsqq}
[\overline{M}_{g+h}(\eta,\beta)]^{vir} \scap c_{2h}(R_{g,h})\scap 
\tau^{-1}(C).
\end{equation}
The latter class (\ref{sqsqsqq}) is an integral
over the virtual class of the fiber
$\tau^{-1}(C)= \overline{M}_{g+h}(C, [C])$. We find:
\begin{equation}
\label{vbn}
C_g(h,X,\beta)= \int_{[\overline{M}_{g+h}(C, [C])]^{vir}}
c_{2h}(R_{g,h}).
\end{equation}
This integral is identical to (\ref{fortt})
except for the different normal bundles $N$ occurring
in the definition of $R_{g,h}$.

The method in Section 2 to compute (\ref{fortt}) also
yields a computation of (\ref{vbn}). As remarked after
equation (\ref{hkl}), the bundle $N$ affects the integral
(\ref{vbn}) through $\int_C c_1(N^*)$:  
\begin{equation}
\label{hkll}
 C_g(h,X,\beta)= \sum_{\tau\in P(h)} \frac{(\int_C 
c_1(N^*))^l}{|\text{Aut}_\tau|} 
\prod_{i=1}^l
\int_{\overline{M}_{h_i,1}} 
\lambda_{h_i} \lambda_{h_i-1} (\sum_{j=0}^{h_i-1}(-1)^j
\lambda_j \psi_1^{h_i-1-j}).
\end{equation}
Since $\int_C c_1(N^*)=2-2g+K_X\cdot \beta$, Theorem 3
follows via the series analysis of Section 2.

\section{\bf Theorem 4}
Let $\pi: \overline{M}_{q,1} \rarr \overline{M}_{q}$ be the universal
curve (for $q\geq 2$). The class $\psi_1$ is the Chern class
of the cotangent line bundle on $\overline{M}_{q,1}$.
The kappa classes are defined by
$\kappa_j= \pi_*(\psi_1^{j+1})$.
Define
$$\beta_{q-2} = \pi_* (\sum_{j=0}^{q-1} (-1)^j \lambda_j \psi_1^{q-1-j})
= \sum_{j=0}^{q-2} (-1)^j \lambda_j \kappa_{q-2-j}.$$
In the notation of Section 2.3, we see:
$$Q(t) = t^2/24 + \sum_{q\geq2} t^{2q} \int_{\overline{M}_{q}}
\lambda_q \lambda_{q-1} \scap \beta_{q-2}.$$
The results of Section 2.3 applied in case $g=0$ prove:
$$\text{exp}(2Q(t)) = \Big( \frac{t/2}{\sin(t/2)}  \Big)^2.$$
After taking the logarithm, we find:
\begin{equation}
\label{fbg}
Q(t)=
\text{log} \Big( \frac{t/2}{\sin(t/2)}  \Big).
\end{equation}
The right series in (\ref{fbg}) may be expanded as
$$
\log \left( \frac{t/2}{\sin(t/2)} \right) =
\sum_{q\ge1} \frac{|B_{2q}|}{(2q)(2q)!}t^{2q} 
$$
by Lemma 3 of [FP].
Faber has computed 
$$\int_{\overline{M}_q} \lambda_q \lambda_{q-1} \scap
\kappa_{q-2} = \frac{1}{2^{2q-1}(2q-1)!!} \frac{|B_{2q}|}{2q}$$
from Witten's conjectures/ Kontsevich's theorem [F2].
It is known $\mathcal{R}^{q-2}(M_q)$ is
exactly 1 dimensional ([F2], [L]). Since
$\lambda_q\lambda_{q-1}$ vanishes when restricted to $\partial 
\overline{M}_q$, we find
$$\beta_{q-2} = \frac{\int_{\overline{M}_{q}}
\lambda_q \lambda_{q-1} \scap \beta_{q-2}}
{\int_{\overline{M}_q} \lambda_q \lambda_{q-1} \scap
\kappa_{q-2}} \cdot \kappa_{q-2}.$$
Theorem 4 now follows from the computation:
$$
\frac{\int_{\overline{M}_{q}}
\lambda_q \lambda_{q-1} \scap \beta_{q-2}}
{\int_{\overline{M}_q} \lambda_q \lambda_{q-1} \scap
\kappa_{q-2}} = \frac{2^{q-1}}{q!}.$$

\vspace{+10 pt}
\noindent
Department of Mathematics \\
\noindent California Institute of Technology \\
\noindent Pasadena, CA 91125 \\
\noindent rahulp@cco.caltech.edu
\end{document}